\documentclass[12pt]{iopart}

\usepackage{iopams}  
\usepackage{amssymb}
\usepackage{amsmath}
\usepackage{amsrefs}
\usepackage [utf8]{inputenc}
\usepackage[T1]{fontenc}
\usepackage{fancybox}
\usepackage{fancyhdr}

\newtheorem{theo}{Theorem}[section]
\newtheorem{lem}[theo]{Lemma}
\newtheorem{theoa}{Theorem}
\newtheorem{prop}[theo]{Proposition}

\newtheorem{cora}{Corollary}[theoa]

\newtheorem{example}[theo]{Example}

\newtheorem{remark}[theo]{Remark}

\newcommand{\NN}{\mathbb{N}}
\newcommand{\PP}{\mathbb{P}}
\newcommand{\RR}{\mathbb{R}}

\newcommand{\ZZ}{\mathbb{Z}}

\newcommand{\CC}{\mathbb{C}}
\newcommand{\TT}{{\mathbb{T}}}
\newcommand{\DD}{{\mathbb{D}}}

\newcommand{\Fb}{\mathbf{F}}
\newcommand{\Cr}{\mathcal{C}}
\newcommand{\Br}{\mathcal{B}}
\newcommand{\Sr}{\mathcal{S}}

\newcommand{\Rep}{\mathrm{Rep}}
\newcommand{\Aut}{\mathrm{Aut}}

\newcommand{\Hom}{\mathrm{Hom}}

\newcommand{\qed}{\begin{flushright} $\square$ \end{flushright}}
\newcommand{\dr}{\partial}

\fancypagestyle{monstyle}{%
	\fancyhf{} 
	                  
	\lhead{Dynamical Green functions and Schrödinger operators}  
 	\rhead{\thepage}
 	\rfoot{{\scriptsize\textcopyright 2013 IOP Publishing Ltd}}
}
\begin{document}

\noindent\shadowbox{\begin{minipage}{\textwidth}\noindent \textbf{Disclaimer :} \textit{this is an author--created, un--copyedited version of an article accepted for publication in} Nonlinearity\textit{. IOP Publishing Ltd is not responsible for any errors or omissions in this version of the manuscript or any version derived from it. The Version of Record is available online at}  doi:10.1088/0951-7715/27/3/527 \textit{.}\end{minipage}}

\title[Dynamical Green functions and Schrödinger operators]{Dynamical Green functions and discrete Schrödinger operators with potentials generated by primitive invertible substitution}

\author{Arnaud Girand}
\address{IRMAR (UMR 6625 du CNRS), Universit\'e de Rennes 1, campus de Beaulieu. B\^at. 22-23. 35042 Rennes Cedex. France.}
\ead{arnaud.girand@univ-rennes1.fr}

\begin{abstract}
In this paper, we set up a "dictionary" between discrete Schrödinger operators and holomorphic dynamics on certain affine cubic surfaces, building on previous work by Cantat, Damanik and Gorodetski. To achieve this, we make use of potential theory: a detailed description of the dynamical Green functions is obtained; then basic results concerning the equilibrium measures and the Green functions of compact subsets of $\CC$ are used to transfer statements from the dynamical context to the Schrödinger one. This provides a new viewpoint on several recent theorems.
\textbf{$ $\\
$ $\\
Latest revision :} \today.
\end{abstract}
%\submitto{\NL}
$ $\\
$ $ \\
\begin{flushright}
{\scriptsize\textcopyright 2013 IOP Publishing Ltd}
\end{flushright}

\maketitle

\section{Set--up and Main Results}

The goal of this paper is to add on previous work by Cantat \cite{BHPS}, Damanik and Gorodetski \cite{Dam4, Dam5} (see also \cite{Cas,S1} for instance) to establish a correspondence between the study of certain discrete Schrödinger operators and the holomorphic dynamics of automorphisms on a family of affine cubic surfaces.

\subsection{Discrete Schrödinger Operators}

\subsubsection{Left shift dynamics.}\label{LSD} Consider the free group on two generators $\Fb_2:= \langle a,b \, | \, \emptyset \rangle$ and let $\varphi \in \Aut(\Fb_2)$ be a positive automorphism, \emph{i.e} such that the images $\varphi(a)$ and $\varphi(b)$ are words in $a$ and $b$ (and thus do not involve the inverse elements $a^{-1}$ and $b^{-1}$). 

Using the action of $\Aut(\Fb_2)$ on the abelianized group $\mathrm{Ab}(\Fb_2) = \ZZ^2$, one can associate a matrix $M_\varphi \in GL_2(\ZZ)$ to $\varphi$. Assume $M_\varphi$ to be hyperbolic, \emph{i.e}: 
\begin{itemize}
\item either $\det(M_\varphi) = 1$ and $\tr(M_\varphi) > 2$;
\item or $\det(M_\varphi) = -1$ and $\tr(M_\varphi)\neq 0$.
\end{itemize}
By replacing $\varphi$ with $\varphi^2:= \varphi \circ \varphi$ ,which is still positive, we can restrict ourselves to the first case; this means that the spectrum of $M_\varphi$ is of the form $\lbrace \lambda, \lambda^{-1} \rbrace$ where $\lambda$ denotes a quadratic integer greater than one.

Let $\Omega$ be the set of finite words on the generators $a$ and $b$, endowed with the topology pertaining to the following distance: 
\begin{displaymath}
d: (u,v) \mapsto \dfrac{1}{\inf \lbrace |n| \, | \, u_n \neq v_n \rbrace + 1}.
\end{displaymath}
The initial automorphism $\varphi$ extends to a substitution $\iota$ over the letters $a$ and $b$ which has an unique "positively infinite" invariant word $u_+ \in \lbrace a , b \rbrace^\NN$.

\begin{example} Let $\zeta$ be the Fibonacci substitution, given by $a \mapsto ab$ and  $b \mapsto a$; its associated matrix $M_\zeta$ is given by: 
\begin{displaymath}
 M_\zeta:= \begin{pmatrix}
1 & 1 \\ 
1 & 0
\end{pmatrix} \in GL_2(\ZZ)
\end{displaymath}
and it fixes the infinite word beginning with $abaababaabaababaababaabaababaabaab \ldots$ This example will be continued throughout the first part of this paper.
\end{example}

Now consider the left shift on $\lbrace a,b \rbrace^\ZZ$: 
\begin{align*}
T:  \lbrace a,b \rbrace^\ZZ & \rightarrow \lbrace a,b \rbrace^\ZZ \\
u & \mapsto (u_{n+1})_{n \in \ZZ}
\end{align*}
and let $W$ be the set of all adherent values for the sequence $(T^p u_+)_{p \geq 0}$ (in other words, it is the $\omega$--limit set $W$ of the $T$--orbit of $u_+$); it is a compact subset of $\lbrace a,b \rbrace^\ZZ$. It is well known (see for instance \cite{Dam1}) that there exists an unique $T$--invariant probability measure $\nu$ on the set $W$ and that the left shift $T$ is ergodic with respect to $\nu$ (see \cite[p.58]{RM} for an outlook on uniquely ergodic maps). 

\subsubsection{Discrete Schrödinger operators.}
Given any word $w \in W$, one can define a potential function: 
\begin{align*}
v_w: \ZZ & \rightarrow \lbrace 0 , 1 \rbrace \\
n & \mapsto \left\lbrace \begin{array}{l}
1 \text{ if } w_n = a \\ 
0 \text{ else } 
\end{array} \right. .
\end{align*}
Consider for any fixed $\kappa \in \RR$ and $w \in W$ the following operator, defined on the space $\ell^2(\ZZ)$ of complex--valued square--summable sequences: 
\begin{align*}
H_{\kappa , w}: \ell^2(\ZZ) & \rightarrow \ell^2(\ZZ)\\
\xi & \mapsto (\xi_{n+1} + \xi_{n-1} + \kappa v_w(n) \xi_n)_{n \in \ZZ}.
\end{align*}
Remark that this operator is self--adjoint and $
\|H_{\kappa,w}\| \leq 2 + |\kappa|;$ therefore its spectrum $\Sigma_{\kappa, w}$ is a subset of the real interval $[-2-|\kappa|, 2 + |\kappa|]$. 

Since $H_{\kappa,w}$ is uniquely ergodic, we can apply the following result due to Kotani and Pastur \cite{Pastur}.

\begin{theo}[Kotani -- Pastur]
$ $\\
There exists a compact set $\Sigma_\kappa \subset [-2-|\kappa|, 2+|\kappa|]$ such that $\Sigma_{\kappa,w} = \Sigma_\kappa$ for all $w \in W$.
\end{theo}
We call the set $\Sigma_\kappa$ the \textbf{almost--sure spectrum} of the operator $H_{\kappa,w}$ with respect to the measure $\nu$.

\begin{remark}
If $H_{\kappa,w}$ were ergodic (non--uniquely), one would have $\Sigma_{\kappa,w} = \Sigma_\kappa$ for $\nu$--almost every $w \in W$, hence the colloquial name of "almost--sure spectrum".
\end{remark}

\subsubsection{Density of states. }Let $H^N_{\kappa,w}$ be the restricted operator ${H_{\kappa,w}}$ to the set $\CC^{\lbrace -N, \ldots , N \rbrace}$ with Dirichlet boundary conditions, meaning we only consider sequences $(\xi_n)_{n\in \ZZ}$ with $-N \leq n \leq N$ such that: 
\begin{itemize}
\item $\xi_n = 0$ for $n \leq -N-1$; 
\item $\xi_n = 0$ for $n \geq N+1$.
\end{itemize}
This gives a self--adjoint endomorphism of $\CC^{2N+1}$; as such it has real eigenvalues $\lambda_0^N, \ldots , \lambda_{2N}^N$. Define the following probability measure: 
\begin{displaymath}
\mu_N^{\kappa,w}:= \dfrac{1}{2N+1} \sum_{j=0}^{2N} \delta_{\lambda_j^N}.
\end{displaymath}

\begin{theo}[Avron -- Simon \cite{AS}]\label{thm:Avron-Simon}$ $
\begin{enumerate}
\item[(i)] For $\nu$--almost every $w \in W$ the sequence $(\mu_N^{\kappa,w})_N$ weakly converges to a probability measure $dk_\kappa$ on $\CC$, called \textbf{density of states};
\item[(ii)] for any continuous function $g: \CC \rightarrow \CC$: 
\begin{displaymath}
\int_\CC g(E) \rmd k_\kappa(E) = \int_{w \in W} \langle g(H_{\kappa,w}) \cdot \delta_0 \, | \, \delta_0 \rangle \rmd \nu(w) \;;
\end{displaymath}
\item[(iii)] the support $\mathrm{supp}(dk_{\kappa})$ of $dk_\kappa$ satisfies $\mathrm{supp}(dk_{\kappa}) = \Sigma_\kappa$.
\end{enumerate}
\end{theo}

\begin{remark}
It is standard to then define the \textbf{integrated density of states} as the distribution function of the probability measure $dk_\kappa$: 
\begin{displaymath}
k_\kappa: E  \mapsto \int_{-\infty}^E \rmd k_\kappa.
\end{displaymath}
\end{remark}

\subsubsection{Lyapunov exponent. }
A hypothetical eigenvalue--eigenvector pair $(E,\xi)$ for $H_{\kappa,w}$ should satisfy the equation: 
\begin{equation}
\forall n \in \ZZ , \quad \xi_{n+1} + \xi_{n-1} + \kappa v_w(n) \xi_n = E \xi_n \: , 
\end{equation}
that is: 
\begin{equation} \label{eigen-eq}
\forall n \in \ZZ, \quad 
\begin{pmatrix}
\xi_{n+1} \\ 
\xi_n
\end{pmatrix} = M_{n,\kappa,w}^E \begin{pmatrix}
\xi_n \\ 
\xi_{n-1}
\end{pmatrix}
\end{equation}
where: 
\begin{displaymath}
 M_{n,\kappa,w}^E:= \begin{pmatrix}
E-\kappa v_w(n) & -1 \\ 
1 & 0
\end{pmatrix} \in SL_2(\CC),
\end{displaymath}
\emph{i.e}  $M_{n,\kappa,w}^E$ is equal to one of the two matrices: 
\begin{displaymath}
M_{\kappa}^E(a):= \begin{pmatrix}
E-\kappa & -1 \\ 
1 & 0
\end{pmatrix} , \quad M_{\kappa}^E(b):= \begin{pmatrix}
E & -1 \\ 
1 & 0
\end{pmatrix}.
\end{displaymath}

Consider the \textbf{Lyapunov exponent}: 
\begin{displaymath}
\gamma_\kappa (E):= \limsup_{N \to \infty} \dfrac{1}{N} \int_W \log \left\| \prod_{n=0}^N M_{n,\kappa,w}^E \right\| \rmd \nu(w).
\end{displaymath}
By Oseledets Theorem, this quantity is well defined and
\begin{displaymath}
\limsup_{N \to \infty} \log \left\| \prod_{n=0}^N M_{n,\kappa,w}^E \right\|
\end{displaymath}
is $\nu$--almost surely constant equal to $\gamma_\kappa(E)$. 

\begin{theo}[see\cite{Dam1} and \cite{CarLa}]\label{IKP}
$ $\\
The Lyapunov exponent is a non--negative function such that : 

\begin{enumerate}
\item[$(i)$] \begin{displaymath}
\gamma_\kappa(E) = \int_{\Sigma_\kappa} \log |E-E^\prime| \rmd k_\kappa(E^\prime);
\end{displaymath}
\item[$(ii)$] \begin{equation}\label{lelong_gamma}
d d^c \gamma_\kappa = 2 \pi dk_\kappa;
\end{equation}
\item[$(iii)$] the almost--sure spectrum satisfies $\Sigma_\kappa = \lbrace \gamma_\kappa = 0 \rbrace$.
\end{enumerate}
\end{theo}
\textit{Proof.}
The first item is the \textbf{Thouless formula} (see \cite[p.340]{CarLa}) and thus, since $\rmd \rmd ^c \log |z-z_0| = 2\pi \delta_{z_0}$, one obtains property $(ii)$. The third result is a theorem due to Ishii, Kotani and Pastur (see \cite{Dam1} for an overview).
\qed

\subsubsection{Green function for the almost--sure spectrum. }First, recall the following definition. Let $U$ be an open set in $\CC$ such that its complement $\CC \setminus U$ is a compact set. A function $g_U: U \rightarrow (0 ,\infty)$ is a \textbf{Green function} for the domain $U$ (alternatively, for the compact $\CC \setminus U$) if:
\begin{enumerate}
\item[$(G1)$] $g_U$ is harmonic;
\item[$(G2)$] the following limit exists: 
\begin{displaymath}
\lim_{z \to \infty} (g_U(z) - \log|z|);
\end{displaymath}
\item[$(G3)$] for all $\xi \in \dr U$, one has: 
\begin{displaymath}
\lim_{z \to \xi} g_U(z) = 0.
\end{displaymath}
\end{enumerate}

\begin{remark}$ $
\begin{enumerate}
\item If $U$ is such an open subset of $\CC$ then its Green function, if it exists, is unique (see \cite[p.182]{Klim}). Moreover, one can replace $(G2)$ with $g_U(z) - \log|z| = \Or(1)$ at infinity.
\item If $U$ has a Green function, there exists a positive real number $C$ such that
\begin{displaymath}
g_U(z) = \log|z| - \log(C) + o(1) \text{ as $z$ goes to infinity}.
\end{displaymath}
The number $C$ is called the \textbf{capacity} of the compact set $\CC \setminus U$. For more details on set capacities, see \cite[p.132]{Ran}.
\item The measure $\rmd \rmd ^c g_U$ is called the \textbf{equilibrium measure} of the compact set $\CC \setminus U$.
\end{enumerate}
 \end{remark}

Consider the open set $U:= \CC \setminus \Sigma_\kappa$; it satisfies $\dr U = \Sigma_\kappa$. We then have the following result, which is well known to experts (see for example \cite[p.979]{Dam4}, remark (g)).

\begin{prop}\label{rem_cap}
$ $
\begin{enumerate}
\item[$(i)$]The Lyapunov exponent $\gamma_\kappa$ is the Green's function for the domain $U$.
\item[$(ii)$] The density of states is the equilibrium measure of $\Sigma_\kappa$.
\item[$(iii)$] The capacity  $\mathrm{Cap}(\Sigma_\kappa)$ of the almost--sure spectrum is one.
\end{enumerate}
\end{prop}
\textit{Proof.} The Thouless formula shows that $\gamma_\kappa: U \rightarrow (0, \infty)$ satisfies condition $(G1)$; moreover, for $E \in \CC$: 
\begin{align*}
\gamma_\kappa (E) - \log|E| & = \int_{\Sigma_\kappa} \log|E-E^\prime| \rmd k_\kappa(E^\prime) - \log|E|\\
& = \int_{\Sigma_\kappa} \log \left| 1 - \dfrac{E^\prime}{E} \right| \rmd k_\kappa(E^\prime)\\
& \xrightarrow[E \to \infty]{}0,
\end{align*}
where the final line follows from the preceding because the function $\log|1-E^\prime / E |$ converges uniformly towards zero on the compact support $\Sigma_\kappa$ of $\rmd k_\kappa$. Therefore condition $(G2)$ holds. Finally, one checks $(G3)$ using Theorem \ref{IKP}. Thus $(i)$ and $(ii)$ hold, using \ref{lelong_gamma} and since $\gamma_\kappa (E) - \log|E| \xrightarrow[E \to \infty]{}0$, one immediately gets $(iii)$.
\qed

\subsection{Holomorphic Dynamics}

\subsubsection{Character variety of the free group on two generators}\label{charac} Let us fix a generating set $\lbrace a ,b \rbrace$ of the free group $\Fb_2$ and consider the algebraic quotient $\chi(\Fb_2)$ of: 
\begin{displaymath}
\Rep(\Fb_2):= \Hom(\Fb_2, SL_2(\CC)) \cong SL_2(\CC) \times SL_2(\CC)
\end{displaymath}
under $SL_2(\CC)$--conjugacy. The variety $\chi(\Fb_2)$ is isomorphic to $\CC^3$ with the following projection map: 
\begin{align*}
\chi: \Rep(\Fb_2) & \mapsto \CC^3\\
\rho & \mapsto (x,y,z) = ( \tr(\rho(a)), \tr(\rho(b)), \tr(\rho(ab))).
\end{align*}
Moreover, if one enforces the condition $\tr([\rho(a),\rho(b)]) = D - 2 \in \CC$ one obtains an affine cubic surface $\Sr_D$, the equation of which is (see \cite{BHPS, CL} for details): 
\begin{displaymath}
x^2 + y^2 + z^2 = xyz +  D.
\end{displaymath}
Let $\varphi$ be an element of $\Aut(\Fb_2)$; then the following defines an automorphism of the surface $\Sr_D$: 
\begin{displaymath}
 f: \chi(\rho) \mapsto \chi(\rho \circ \varphi^{-1}).
\end{displaymath} 
Since the group  $\Aut(\Fb_2)$ acts on $\mathrm{Ab}(\Fb_2) = \ZZ^2$ one can set  
\begin{displaymath}
M_f = \begin{pmatrix}
p & q \\ 
r & s
\end{pmatrix} \in GL_2(\ZZ)
\end{displaymath}
to be the matrix corresponding to $\varphi^{-1}$ and if $A:= \rho(a)$,  $B:= \rho(b)$ for some $\rho \in \Rep(\TT^2_1)$ then: 
\begin{equation} \label{egalite_PGL2}
f ( \chi(\rho) ) = ( ( \tr(A^p B^q), \tr(A^r B^s), \tr(A^pB^qA^rB^s)).
\end{equation}
This gives us an action of $GL_2(\ZZ)$ on $\Sr_D$ whose kernel contains $\pm I_2$; therefore $PGL_2(\ZZ)$ acts on the surface $\Sr_D$. 

Using (\ref{egalite_PGL2}) and Fricke--Klein's formulas, one sees that $f$ is a polynomial automorphism of $\Sr_D$; in the following, we will denote by $\Br$ the subgroup of $\Aut(\Sr_D)$ formed by such mappings $f$. We will say that an automorphism $f \in \Br$ is hyperbolic if one of the next two conditions holds: 
\begin{itemize}
\item either $\det(M_f) = 1$ and $\tr(M_f) > 2$;
\item or $\det(M_f) = -1$ and $\tr(M_f)\neq 0$.
\end{itemize}

\begin{example}
For the Fibonacci substitution $\zeta$, consider the automorphism $f$ associated with $\zeta
$. Then: 
\begin{displaymath}
M_f = M_{\zeta^{-1}} = \begin{pmatrix}
0 & 1 \\ 
1 & -1
\end{pmatrix}.
\end{displaymath}
Since $\det(M_f) = -1$ and $\tr(M_f) = -1 \neq 0$ the morphism $f$ is in fact hyperbolic. It is given by:
\begin{displaymath}
f(x,y,z) = (y,xy-z,x).
\end{displaymath}
\end{example}
Denote by $\overline{\Sr_{D}} \subset \PP^3$ the compactified surface: 
\begin{displaymath}
w(x^2 + y^2 + z^2) = xyz + w^3D,
\end{displaymath}
where $[x:y:z:w]$ are homogeneous coordinates on the projective space $\PP^3$. Its intersection with the plane at infinity $\lbrace w = 0 \rbrace$ is equal to the "triangle at infinity" $\Delta = \lbrace xyz = 0 \rbrace$. Thus $\Aut(\Sr_{D})$ embeds into the group of birational transformations of $\overline{\Sr_{D}}$. The dynamics at infinity of the hyperbolic elements in $\Br$ is quite rich, as we will see throughout this paper; first, we have the following result.

\begin{prop}[see \cite{BHPS,CL,EH}]
$ $\\
Let $f \in \Br$ be a hyperbolic automorphism. Then $f$ extends to a birational transformation of $\overline{\Sr_D}$ and: 
\begin{enumerate}
\item[(i)] $f$ has an unique indeterminacy point $v_{-}$ which is either $[1:0:0:0]$, $[0:1:0:0]$ or $[0:0:1:0]$;
\item[(ii)] the mapping $f$ contracts $\Delta \setminus \lbrace v_- \rbrace$ onto the indeterminacy point $v_+$ of $f^{-1}$;
\item[(iii)] up to conjugacy by an element of $\Br$, one can assume $v_+$ to be distinct from $v_-$.
\end{enumerate}
\end{prop}

\begin{remark}\label{rk_EH}
Èl'Huti \cite{EH} gave a detailed description of the automorphism group $\Aut(\Sr_D)$; in particular, he proved that $\Br$ has finite index in $\Aut(\Sr_D)$.
\end{remark}

\subsubsection{Main theorem on dynamical Green functions. } Fix a hyperbolic automorphism $f \in \Br$ for which $v_+ \neq v_-$ and denote by $\lambda$ the spectral radius of $M_f$. We now try to understand the escape rate at infinity in the unbounded orbits under $f$. First, a theorem by Dloussky \cite{Dlou} combined with work by Cantat \cite{BHPS} (see also \cite{Fav}) yields the following result, which will be essential to our study of the dynamics of $f$ at infinity.

\begin{prop}[\cite{BHPS}, theorem 3.1 p. 423]\label{prop_dloussky}
$ $\\
There exists a matrix $N_f \in GL_2(\ZZ)$ with non--negative entries which is conjugate to $M_f$ in $PGL_2(\ZZ)$, an open neighbourhood $U$ of $v_+$ in $\overline{\Sr_{D}}$ and a biholomorphism $\psi_f^+: \DD \times \DD \rightarrow U$ such that: 
\begin{enumerate}
\item[(i)] $\psi_f^+(0,0) = v_+$;
\item[(ii)] for all $(u,v) \in \DD^* \times \DD^*$ one has: 
\begin{displaymath}
\psi_f^+((u,v)^{N_f}) = f ( \psi_f^+(u,v)),
\end{displaymath}
where $(u,v)^{N_f}$ denotes the monomial action of $N_f$ on the pair $(u,v)$, \emph{i.e} if 
\begin{displaymath}
N_f = \begin{pmatrix}
p & q \\ 
r & s
\end{pmatrix} 
\end{displaymath}
then $(u,v)^{N_f} = (u^p v^q, u^r v^s)$.
\end{enumerate}
As a consequence, if $m \in \Sr_{D}$ has unbounded forward orbit under $f$, then $f^n(m)$ goes to $v_+$ at infinity.
\end{prop}

Before stating our main result regarding dynamical Green functions, let us set a few conventions: 
\begin{itemize}
\item define the \textbf{filled Julia set} $K^+(f)$ as follows: 
\begin{displaymath}
K^+(f):= \lbrace m \in \Sr_{D} \, | \, \exists M > 0, \, \forall n \geq 0, \, \|f^n(m)\| \leq M \rbrace \; ,
\end{displaymath}
where $\|.\|$ denotes the standard euclidean norm on $\CC^3$; 
\item set $\alpha, \beta \in \RR_+^*$ to be the coordinates of the projection of the vector $ \begin{pmatrix}
1 \\ 
1
\end{pmatrix} $ on the eigenline for $N_f$ associated with the maximal eigenvalue of $M_f$ (and so of $N_f$).
\end{itemize}

\begin{theoa}[Dynamical Green function]\label{lem_equiv}
$ $\\
Let $f \in \Br$ be a hyperbolic element and let $m \in \Sr_D$. Then the following quantity is well defined: 
\begin{displaymath}
G_f^+: m \mapsto \lim_{n \to \infty} \dfrac{1}{\lambda^n}\log^+ \|f^n(m)\| \, ,
\end{displaymath}
and: 
\begin{enumerate}
\item[(i)] the function $G_f^+$ is pluriharmonic (resp. plurisubharmonic) on the complement of the filled Julia set $K^+(f)$ in $\Sr_{D}$ (resp. on $\Sr_D$) and takes non--negative values;
\item[(ii)] the zero set of $G_f^+$ is  $K^+(f)$;
\item[(iii)] the following relation holds: 
\begin{equation}
G_f^+ \circ f = \lambda G_f^+;
\end{equation}
\item[(iv)] if $m = \psi_f^+(u,v) \in \psi_f^+(\DD^* \times \DD^*)$, then: 
\begin{equation}\label{equiv2}
G_f^{+}(m) = -\alpha \log|u| - \beta  \log|v| 
\end{equation}
\item[(v)] the function $G_f^+$ is locally Hölder--continuous.
\end{enumerate}
\end{theoa}

\begin{example}
In the Fibonacci case, $M_f$ is conjugate to
\begin{displaymath}
N_f= \begin{pmatrix}
1& 1 \\ 
1 & 0
\end{pmatrix} \text{ in } PGL_2(\ZZ).
\end{displaymath}
The eigenvalues of $M_f$ (and so of $N_f$) are $$ \phi:= \dfrac{1+\sqrt{5}}{2} \quad \text{ and } \quad \overline{\phi}:= \dfrac{1-\sqrt{5}}{2}.$$ and the corresponding eigenlines for $N_f$ are spanned by $\begin{pmatrix}
\phi  \\ 
1 
\end{pmatrix}$ and $\begin{pmatrix}
\overline{\phi}  \\ 
1 
\end{pmatrix}$. Thus, since: 
\begin{displaymath}
\begin{pmatrix}
1  \\ 
1 
\end{pmatrix}=\dfrac{1-\overline{\phi}}{\sqrt{5}}\begin{pmatrix}
\phi  \\ 
1 
\end{pmatrix}+ \dfrac{\phi-1}{\sqrt{5}}
\begin{pmatrix}
\overline{\phi}  \\ 
1 
\end{pmatrix}\, , \end{displaymath}
one has $\alpha =\dfrac{1-\overline{\phi}}{\sqrt{5}} \phi = \dfrac{\phi-1}{\sqrt{5}} $ and $\beta = \dfrac{1-\overline{\phi}}{\sqrt{5}} $. Moreover, we have in this case $v_+ = [0:1:0:0]$.
\end{example}

\subsection{Applications to Discrete Schrödinger Operators}

\subsubsection{Schrödinger curve} Consider the following cubic surface in $\CC^3$, for some fixed $\kappa \in \RR$: 
\begin{displaymath}
(\Sr_{4+\kappa^2}) \qquad x^2 + y^2 + z^2 = xyz + 4 + \kappa^2 \quad; 
\end{displaymath}
this is a connected smooth ($if \kappa \neq 0$) affine surface, containing what we call its Schrödinger curve: \begin{align*}
s: \CC & \rightarrow \Sr_{4+\kappa^2}\\
E & \mapsto (E-\kappa, E , E(E-\kappa) - 2).
\end{align*}

\begin{remark}
The function $s$ is in fact the trace map associated with the matrices $M_{n,\kappa,w}^E$. Namely, one has $s(E) = ( \tr(M_{\kappa}^E(a)), \tr(M_{\kappa}^E(b)), tr(M_{\kappa}^E(b)M_{\kappa}^E(a)))$.
\end{remark}

Starting from our automorphism $\varphi \in \Aut(\Fb_2)$ (cf. \ref{LSD}) with associated substitution $\iota$, we obtain a polynomial automorphism $f$ of $\Sr_{4+\kappa^2}$ associated with $\varphi^{-1}$ (cf. \ref{charac}); one can then explicitly compute it using the formula $ f(\chi(\rho)) = \chi(\rho \circ \varphi)$ and so its restriction to the Schrödinger curve is: 
\begin{displaymath}
\forall E \in \CC, \quad f(s(E)) = ( \tr(M_{\kappa}^E (\iota(a))),\tr(M_{\kappa}^E (\iota(b))),\tr(M_{\kappa}^E(\iota (ab)))),
\end{displaymath}
where, if $u = (u_1, \ldots, u_n) \in \lbrace a,b \rbrace^n$, then: 
\begin{displaymath}
M_{\kappa}^E(u):= \prod_{i=0}^{n-1} M_{\kappa}^E (u_{n-i}).
\end{displaymath}
Since $\varphi$ is hyperbolic, $f$ is a hyperbolic automorphism of $\Sr_{4+\kappa^2}$. We then have the following result \cite{Dam2} (see also some earlier work by Süt\H{o} \cite{S1,S2}).

\begin{prop}[Damanik \cite{Dam2}, theorem 2.1 p. 399]\label{damanik_th}
$ $\\
If $f$ is the polynomial automorphism of $\Sr_{4+\kappa^2}$ associated with a positive hyperbolic substitution $\iota$ on two letters, then the almost--sure spectrum $\Sigma_\kappa$ satisfies:
\begin{displaymath}
\Sigma_\kappa = s^{-1} (K^+(f)).
\end{displaymath}
\end{prop}

\subsubsection{"Dictionary" Between Holomorphic Dynamics and Schrödinger Operators}
We now move on to our second result. Since the subgroup $\Br$ has finite index in $\Aut(\Sr_{4+\kappa^2})$ (cf. remark \ref{rk_EH}) and $f$ has infinite order we can suppose, up to replacing it with some iterate $f^{n_0}$ that $f \in \Br$; thus, we will be able to exploit theorem $A$ to obtain the following result.

\begin{theoa}\label{theo_B}
$ $\\
Let $\iota$ be a positive hyperbolic substitution over the letters $a$ and $b$ and let $f\in \Br$ be the associated automorphism of $\Sr_{4+\kappa^2}$. Then for $E \in \CC$: 
\begin{displaymath}
\gamma_\kappa (E) = \dfrac{1}{\alpha + \beta} G_f^+ (s(E)),
\end{displaymath}
where $\alpha, \beta \in \RR_+^*$ are the same as in Theorem \ref{lem_equiv}.
\end{theoa}

\begin{remark}
Proposition \ref{damanik_th} was mostly a qualitative one, concerning the boundedness of the orbit alone. Here, using our Theorem A, we get tools to estimate the escape rate at infinity thus obtaining a more quantitative result.
\end{remark}

This, combined with previous work by Cantat, Damanik and Gorodetski, allows us to work out the following "dictionary".
$ $\\
\fbox{\begin{minipage}{0.915\textwidth}
\begin{tabular}{c|c}
\textbf{Discrete Schrödinger operators}  & \textbf{ Holomorphic dynamics on $\Sr_{4+\kappa^2}$} \\ 
\hline & \\
 Almost--sure spectrum $\Sigma_\kappa$ &   Julia set $K^+(f)  $\\ 
 Lyapunov exponent $\gamma_\kappa$ &   Dynamical Green's function $G_f^+  $\\ 
 Density of states $\rmd k_\kappa$ &  Green's current $  T^+_f$  \\ 
Thouless formula & $T^+_f = \rmd \rmd ^c G_f^+$  \\ 
 $\gamma_\kappa$ and $k_\kappa$  Hölder--continuous near $\Sigma_\kappa$ & $G_f^+$ locally Hölder--continuous  \\ 
 Avron and Simon Theorem~\ref{thm:Avron-Simon} & Convergence to $T^+_f$
\end{tabular} 
\end{minipage}}
$ $\\

More precisely, one goes from the right-hand side of this table to the left by taking pull-backs with the Schr\"odinger curve $s\colon \CC\to S_{4+\kappa^2}$; for 
instance, the first line is Damanik's Proposition~\ref{damanik_th}, and the second is our Theorem~B. Similarly, the H\"older
continuity of $\gamma_\kappa$ corresponds to the H\"older continuity of $G_f^+$ (obtained in Theorem~A); we shall see in Section~\ref{par:consequence-holder} that it implies directly H\"older continuity of the integrated density of states. 
The last line of this table is less precise: this is explained in paragraph~\ref{par:consequence-AS}.

It is to be noted that the literature pertaining to the interactions between real and complex dynamics and discrete Schrödinger operators has been quickly expanding these last few years, mainly under the guidance of D. Damanik and A. Gorodetski. See for instance \cite{Mei,Y1,Y2}.

\section{Dynamical Green Functions}

\subsection{Preliminary Computations}
\subsubsection{Geometry of $\Sr_D$ at infinity}\label{Geom}
In order to measure the escape rate at infinity of a point with unbounded orbit, we will now study the behaviour of $\log\|m\|$ when $m=(x,y,z)\in\Sr_{D}$ goes to $v_+$, where $\|\cdot\|$ denotes the euclidean norm on $\CC^3$. For the sake of clarity, suppose (our problem being symmetric with respect to $x,y$ and $z$) that $v_+$ is the point $[0:0:1:0]$; in a neighbourhood of $v_+$, $\Sr_{D}$ can be seen, using the chart $\lbrace z \neq 0 \rbrace$, as the surface: 
\begin{equation}
(X^2 + Y^2 + 1) W = XY + DW^3,
\end{equation}
where $X:= x/z$, $Y:= y/z$ and $W:= w/z$. Equivalently, this can be written as follows: 
\begin{equation}
W = XY + DW^3 + W^2(AX+BY+C) + W(X^2 + Y^2).
\end{equation}
Using these new coordinates $(X,Y,W)$, $v_+$ corresponds to the point at origin $(0,0,0)$ and one has: 
\begin{align*}
\log \|m\| =& \dfrac{1}{2} \log \left( \left| \dfrac{X}{W} \right|^2 + \left| \dfrac{Y}{W} \right|^2 + \dfrac{1}{|W|^2} \right)\\
 = &- \dfrac{1}{2} \log ( |W|^2 ) + \dfrac{1}{2} \log( |X|^2 + |Y|^2 + 1 )\\
 = &- \dfrac{1}{2} \log ( |XY + DW^3 + W^2(AX+BY+C) + W(X^2 + Y^2)|^2 ) \\
& + \dfrac{1}{2} \log( |X|^2 + |Y|^2 + 1 ).
\end{align*}
Using Taylor's approximation one gets: 
\begin{displaymath}
\log \|m\| = -  \log ( |XY| ) + g(X,Y,W),
\end{displaymath}
where $g$ is bounded in a neighbourhood of $(0,0,0)$. Now, for $m$ close enough to $v_+$ one can apply the biholomorphism $\psi_f^+$ to get $(u,v):={\psi_f^+}^{-1}(m)$ and use the following lemma.

\begin{lem}\label{lem_1}
$ $\\
There exists a germ of bounded function $h$ such that for all $(u,v) \in \DD^* \times \DD^*$: 
\begin{displaymath}
\log\|\psi_f^+(u,v)\| = - \log|uv| + h(u,v).
\end{displaymath}
\end{lem}
\textit{Proof.}Using Taylor's theorem at the origin one gets: 
\begin{displaymath}
\psi_f^+(u,v) = v_+ + L(u,v) + R(u,v),
\end{displaymath}
where $L$ is the linear part of $\psi_f^+$ at the origin and $R$ is a smooth bounded function on $\DD \times \DD$ such that $R(u,v) = \Or(\|(u,v)\|^2)$. Since $\psi_f^+$ is a conjugacy between the dynamics of $f$ and $N_f$ and since $f$ (resp. $N_f$) only contracts the axes $\lbrace X=0 \rbrace$ and $\lbrace Y=0 \rbrace$ (resp. $\lbrace u=0 \rbrace$ and $ \lbrace v=0 \rbrace$) on the origin then $L=d\psi_f^+(0,0)$ must be of the form
\begin{displaymath}
\begin{pmatrix}
r_1 & 0 \\ 
0 & r_2
\end{pmatrix} \text{ or } \begin{pmatrix}
0 & r_1 \\ 
r_2 & 0
\end{pmatrix}.
\end{displaymath}
Therefore, there exists a bounded function $h$ on $\DD \times \DD$ such that: 
\begin{equation}\label{egalite_lisse}
\log \| \psi_f^+(u,v) \| = -  \log (|uv| ) + h(u,v),
\end{equation}
\qed

\subsubsection{Estimate at infinity}\label{MDI}
Since $M_f$ is hyperbolic, one can assume (replacing $M_f$ with $M_{f^2} = M_f^2$) that it has eigenvalues $\lambda$ and $\lambda^{-1}$, with $\lambda$ a real number greater than one. Now consider the following quantity, for $n \geq 0$ and $m$ with unbounded forward orbit, chosen sufficiently close to $v_+$ (\emph{i.e} in ${\psi_f^+}^{-1}(\DD \times \DD)$): 
\begin{displaymath}
\dfrac{1}{\lambda^n} \log \|f^n(m) \|.
\end{displaymath}
Let $(u_n,v_n):= (u,v)^{N_f^n}$; using the previous lemma one gets:
\begin{displaymath}
\dfrac{1}{\lambda^n}\log \| f^n(m) \| = - \dfrac{1}{ \lambda^n } \log (|u_nv_n| ) + \dfrac{1}{\lambda^n} h(u_n	,v_n).
\end{displaymath}
Since $ \dfrac{1}{\lambda^n} h(u_n	,v_n) \xrightarrow[n \to \infty]{}	0$, we want to understand the behaviour at infinity of the following quantity: 
\begin{displaymath}
- \dfrac{1}{ \lambda^n } \log (|u_nv_n| ).
\end{displaymath}

\begin{lem}\label{lem_2}
$ $\\
The following estimate holds, as $n$ goes to infinity: 
\begin{equation}\label{equiv1}
\dfrac{1}{\lambda^n}\log |u_n v_n|\xrightarrow[n \to \infty]{} (\alpha \log |u| + \beta \log |v|),
\end{equation}
where $\alpha, \beta \in \RR_+^*$ are the coordinates of the projection of the vector $ \begin{pmatrix}
1 \\ 
1
\end{pmatrix} $ on the eigenline for $N_f$ associated with $\lambda$.
\end{lem}
\textit{Proof.}Since $(u,v) \in \DD^* \times \DD^*$ one can set $(\rme^s,\rme^t):= (u,v)$ with: $$s,t \in  \lbrace z \in \CC \, | \, \Re(z) < 0 , \Im(z) \in (-\pi,\pi] \rbrace.$$ Then it is just a matter of describing the behaviour of $|uv| = |\rme^{s+t} | = \rme^{\Re(s+t)}$ under $N_f$, which acts linearly on the coordinates $(s,t)$. A computation thus yields: 
\begin{displaymath}
\dfrac{1}{\lambda^n }\log |u_n v_n| \xrightarrow[n \to \infty]{}   (\alpha \Re(s) + \beta \Re(t)) =  (\alpha \log|u| + \beta\log|v|).
\end{displaymath}
\qed

Using lemmas \ref{lem_1} and \ref{lem_2}, one gets the following estimate: 
\begin{equation}\label{equiv_final}
\dfrac{1}{\lambda^n} \log\|f^n(m)\| \xrightarrow[n \to \infty]{} -  (  \alpha \log |u| + \beta \log |v|).
\end{equation}
Note that this only holds for $m$ sufficiently near $v_+$, \emph{i.e} for $(u,v)$ in $\DD^* \times \DD^*$.

\subsection{Proof of Theorem A} 
First remark that if $m \in K^+(f)$ then  it is clear that $G_f^+(m)$ is well defined and equal to $0$.

Now consider $m \notin K^+(f)$; up to replacing $m$ with some $f^{n_0}(m)$, one can assume that $m$ is sufficiently near $v_+$ so that one can set $(u,v):= {\psi_f^+}^{-1} (m)$ and $(u_n, v_n):= (u,v)^{N_f^n}$. Since $f^n(m)\xrightarrow[n \to \infty]{}v_+$, for $n$ large enough, $\log^+\|f^n(m)\| = \log\|f^n(m)\|$. Applying the estimate (\ref{equiv_final}) then yields, using the same notations as before: 
\begin{align*}
 \dfrac{1}{\lambda^n} \log\|f^n(m)\| & \xrightarrow[n \to \infty]{} -  (  \alpha \log |u| + \beta \log |v|) 
\end{align*}
We have thus proved that $G_f^+$ is well defined and that $(ii)$ and $(iv)$ hold. Moreover, the estimate (\ref{equiv_final}) implies that: 
\begin{equation}\label{egalite_exacte}
\forall (u,v) \in \DD^* \times \DD^*, \qquad G_f^+ \circ \psi_f^+(u,v) = - \alpha \log |u| - \beta \log |v|.
\end{equation}

$(i)$ Let $H$ be a compact set in $\Sr_{D}$, $m \in H$ and $n,p \geq 0$. If $ m \in K^+(f)$ then we clearly have uniform boundedness. Else, $f^n(m)\xrightarrow[n \to \infty]{} v_+$ and so for $n$ large enough $ m =\psi_f^+(u,v)$ with $(u,v) \in \DD^* \times \DD^*$ and $\|f^n(m)\| > 1$, therefore $\log^+\|f^n(m)\| = \log\|f^n(m)\|$ and: 
\begin{align*}
\left| \dfrac{1}{\lambda^{n+p}}\log^+ \|f^{n+p}(m)\| -  \dfrac{1}{\lambda^{n}}\log^+ \|f^{n}(m)\|\right| & = \dfrac{1}{\lambda^{n+p}}\left| \log \|f^{n+p}(m)\| -  \lambda^p\log \|f^{n}(m)\|\right| 
\end{align*}
Since we just proved that there exists a constant $C_m = (\alpha \log|u| + \beta \log|v|)$ depending only on the orbit of $m$ such that $\log \|f^n(m)\| = C_m \lambda^n + \lambda^n \varepsilon_m(n) $, with $\varepsilon_m(n) \xrightarrow[n \to \infty]{}0$ so: 
\begin{align*}
\left| \dfrac{1}{\lambda^{n+p}}\log^+\|f^{n+p}(m)\| -  \dfrac{1}{\lambda^{n}}\log^+ \|f^{n}(m)\|\right| & =\dfrac{1}{\lambda^{n+p}}  \left| \varepsilon_m(n+p) - \varepsilon_m(n) \right|. 
\end{align*}
As $|\varepsilon_m (n)| \xrightarrow[n \to \infty]{}0$, then for all positive $\eta$ and $m \in H$, there exists $N_m \in \NN$ such that:  
\begin{displaymath}
\forall n \geq N_m, \quad |\varepsilon_m (n)| \leq |\varepsilon_m (N_m)| < \eta
\end{displaymath}
hence: 
\begin{displaymath}
\left| \dfrac{1}{\lambda^{n+p}}\log^+\|f^{n+p}(m)\| -  \dfrac{1}{\lambda^{n}}\log^+ \|f^{n}(m)\|\right| \leq  2 | \varepsilon_m(N_m) |.
\end{displaymath}
Since $C_m$ and $\log \|f^n(m)\|$ are continuous with respect to $m$ (cf. $(iv)$), $$m \mapsto \varepsilon_m (N_m) = \lambda^{-N_m} ( \log \|f^{N_m}(m)\| - C_m )$$ is continuous. Using the compactness of $H$, there exists $m_0 \in H$ such that: 
\begin{displaymath}
\sup_{m \in H} \varepsilon_m (N_m) = \varepsilon_{m_0} (N_{m_0})
\end{displaymath}
where: 
\begin{align*}
0 \leq  2 | \varepsilon_m(N_m) | \leq 2  | \varepsilon_{m_0}(N_{m_0}) | < 2 \eta.
\end{align*}
The sequence defining $G_f^+$ thus converges uniformly on all compact subsets in $\Sr_{D}$ and so the limit function inherits the pluri(sub)harmonic properties of its terms. \\
$(iii)$ This stems from the fact that if $m \in \Sr$ then: 
\begin{displaymath}
\dfrac{1}{\lambda^n}\log^+\|f^n(f(m))\| = \dfrac{1}{\lambda^n}\log^+\|f^{n+1}(m)\| = \lambda \left(\dfrac{1}{\lambda^{n+1}}\log^+\|f^{n+1}(m)\| \right).
\end{displaymath}
$(v)$ Here we adapt work by Fornaess and Sibony \cite{FS}. Since $G_f^+$ is $\Cr^1$ outside any neighbourhood of $K^+(f)$ it is Hölder--continuous there. Now let $z_1 \in \Sr_D$ and $z_0 \in K^+(f)$ be such that: 
\begin{displaymath}
d(z_1, K^+(f)) = \|z_1 - z_0\|.
\end{displaymath}
If $z_1 \in K^+(f)$, there is nothing to show. Else, note that by definition of the filled Julia set there exists $R_0 > 0$ such that: 
\begin{displaymath}
\forall n \in \NN, \quad \|f^n(z_0) \| \leq R_0.
\end{displaymath}
Let us consider a positive real number $R \geq R_0 + 1$ and set: 
\begin{displaymath}
N:= \min\lbrace n \geq 0 \, | \, \| f^n ( z_1 ) \| > R \rbrace < \infty; 
\end{displaymath}
thus: 
\begin{align*}
| \: \| f^N(z_1) \| - \| f^N (z_0) \| \: | & \leq \| f^N(z_1) - f^N(z_0) \|\\
& \leq \sup_{\|z\| \leq R} \|df(z)\| \| f^{N-1}(z_1) - f^{N-1}(z_0) \| \text{ because } \| f^{N-1}(z_1)\| \leq R\\
& \qquad  \vdots \\
& \leq ( \sup_{\|z\| \leq R} \|df(z)\| )^N\|z_1 - z_0 \|\\
& \leq  ( \sup_{\|z\| \leq R} \|df(z)\| )^N d(z_1, K^+(f)).
\end{align*}
Hence, if one sets: 
\begin{displaymath}
H(R):= \sup_{\|z\| \leq R} \|df(z)\|
\end{displaymath}
one gets:
\begin{align*}
1 \leq R - R_0  & \leq | \: \| f^N(z_1) \| - \| f^N (z_0) \| \: | \\
& \leq H(R)^N d(z_1, K^+(f))
\end{align*}
thus $H(R)^Nd(z_1, K^+(f)) \geq 1$. Setting $\gamma:= \dfrac{\log(\lambda)}{\log(H(R))}$ one has: 
\begin{equation}\label{briend1}
\dfrac{1}{\lambda^N} \leq d(z_1,K)^{\gamma}.
\end{equation}
Using $(iii)$ one gets:
\begin{align*}
G_f^+(z_1) & = \dfrac{1}{\lambda^N} G_f^+ \circ f^N(z_1)\\
& \leq  \dfrac{1}{\lambda^N} \sup_{\|z \| \leq R} G_f^+ \circ f (z) \text{ because } \| f^{N-1}(z_1) \| \leq R \\
& \leq d(z_1, K^+(f))^\gamma \sup_{\|z \| \leq R} G_f^+ \circ f (z) \text{ by (\ref{briend1}).}
\end{align*} 
Let: 
\begin{displaymath}
C:= \sup_{\|z \| \leq R} G_f^+ \circ f (z)
\end{displaymath}
one then has, \emph{in fine}: 
\begin{equation}\label{briend2}
G^+_f(z_1) \leq C d(z_1, K^+(f) )^\gamma
\end{equation}
for any point $z_1 \in \Sr_{D}$.\qed

\begin{remark}
Using the notations of paragraph \ref{MDI}, we can estimate the local coordinates $(X,Y)$ around $v_+$ as follows (up to a permutation of $u$ and $v$ in the linear part): 
\begin{displaymath}
(X,Y) = (r_1 u, r_2 v) + R(u,v).
\end{displaymath}
Therefore, we have, as $m$ goes to $v_+$: 
\begin{equation}\label{equiv_XY}
G_f^+ (m) = - \alpha \log|X| - \beta \log |Y| - \log|r_1^\alpha r_2^\beta | + o(1).
\end{equation}

\end{remark}

\begin{remark} Replacing $f$ with its inverse $f^{-1}$, one can define the negative dynamical Green function: 
\begin{displaymath}
G_f^-:= \lim_{n \to \infty} \dfrac{1}{\lambda^n} \log^+ \|f^{-n}(m) \|.
\end{displaymath}
Our main result extends to this function.
\end{remark}

\subsection{Corollaries}

We can now consider the closed positive current \cite{BHPS} associated with $G_f^+$, namely: 
\begin{displaymath}
T_f^+:= \rmd \rmd ^c G_f^+ = 2\rmi \dr \bar{\dr} G_f^+
\end{displaymath}
which satisfies the following: 
\begin{displaymath}
f^* T_f^+ = \lambda T_f^+
\end{displaymath}
and has support in the Julia set $J^+(f):= \dr K^+(f)$. 
\setcounter{theoa}{1}
\begin{cora}
$ $\\
Let $f \in \Br$ be a hyperbolic element and $m \in \Sr_D$. Then there exists a neighbourhood $U$ of $v_+$ in $\overline{\Sr_{D}}$ such that: 
\begin{displaymath}
{\rmd \rmd ^c G_f^+}_{|_{U}} = - 2 \pi \left( \alpha \int_{X=0} + \beta \int_{Y=0} \right),
\end{displaymath}
where $\alpha, \beta \in \RR_+^*$ and $(v_+,X,Y)$ are the same as in Theorem \ref{lem_equiv}.
\end{cora}
\textit{Proof.}Let $U:= {\psi_f^+}^{-1}(\DD^* \times \DD^*)$; then using (\ref{egalite_exacte}) and (\ref{equiv_XY}) one gets: 
\begin{displaymath}
{\rmd \rmd ^c G_f^+}_{|_{U}} = \rmd \rmd ^c ( - \alpha \log|u| - \beta \log |v| ) = \rmd \rmd ^c ( - \alpha \log|X| - \beta \log |Y| ).
\end{displaymath}
The result then follows from the Lelong--Poincaré lemma. 
\qed

\section{From Holomorphic Dynamics to Schrödinger Operators}

\subsection{Proof of Theorem B}

Consider the function: 
\begin{align*}
g: \CC\setminus  \Sigma_\kappa  & \rightarrow (0, \infty)\\
E & \mapsto G_f^+ \left( s \left( E \right) \right) \quad;
\end{align*}
our aim is to show that it is (up to a multiplicative constant) the Green's function of the domain $U:= \CC\setminus  \Sigma_\kappa $, thus proving the theorem. Since $G_f^+$ is psh, condition $(G1)$ holds and $(G3)$ is a direct consequence of Damanik's result (Proposition \ref{damanik_th}).

Using Fricke--Klein's formulas and relation (\ref{egalite_PGL2}), one shows using induction that $f$ contracts the triangle at infinity $\Delta$ on the point $v_+ = [0:0:1:0]$. Using (\ref{equiv_XY}), one then gets 
\begin{displaymath}
g(E) - \left( - \alpha \log|x| - \beta \log|y| - \log|C|\right) \xrightarrow[E \to \infty]{} 0 
\end{displaymath}
where $C \in \CC$ and $s \left( E \right) = [x:y:1:1]$. One also has: 
\begin{align*}
s([E: t]) = [Et - t^2\kappa: Et: E^2 - Et\kappa -2t^2: t^2],
\end{align*}
hence, using the chart $\lbrace z \neq 0 \rbrace$: 
\begin{align*}
s(E) &= s \left( [E:1] \right) \\
&= \left( \dfrac{E - \kappa }{E^2 - E\kappa -2},\dfrac{E}{E^2 - E\kappa -2}, \dfrac{1}{E^2 - E\kappa -2} \right)\\
&= \left( \dfrac{1}{E }\left(\dfrac{1 - \kappa / E }{1 - \kappa /E  -2 / E^2}\right),\dfrac{1}{E }\left(\dfrac{1 }{1 - \kappa /E  -2 / E^2}\right), \dfrac{1}{E^2 - E\kappa -2} \right)
\end{align*}
Thus the following limit exists: 
\begin{align*}
\lim_{E \to \infty} \left( g(E) -  (\alpha + \beta) \log|E|\right). 
\end{align*} \qed

\begin{remark}
 Using Proposition \ref{rem_cap}, one has: 
 \begin{displaymath}
 \lim_{E \to \infty} \left(g(E) -  (\alpha + \beta) \log|E|\right) = - \log \mathrm{Cap}(\Sigma_\kappa) = 0.
 \end{displaymath}
 \end{remark} 
 
 \subsection{Consequences}
 
 Theorem \ref{theo_B} yields a few interesting corollaries, further detailing the entwining between certain dynamical invariants and discrete Schrödinger operators.

\subsubsection{H\"older continuity (see also \cite{Dam6, Dam7})}\label{par:consequence-holder}

\setcounter{theoa}{2}
\begin{cora}
$ $\\
One has the following results: 
\begin{enumerate}
\item[$(i)$] $s^*( \rmd \rmd ^c G_f^+) =  2\pi(\alpha+\beta) \rmd k_\kappa$;
\item[$(ii)$] the functions $\gamma_\kappa$ and $k_\kappa$ are Hölder--continuous near $\Sigma_\kappa$, with the same H\"older exponent $\tau$;
\item[$(iii)$] the density of states does not charge sets with Hausdorff dimension less than $\tau$. In particular, the Hausdorff dimension of the almost--sure spectrum is strictly positive.
\end{enumerate}
\end{cora}
\textit{Proof.} The first assertion follows from (\ref{lelong_gamma}). 
To prove property $(ii)$, we reproduce an argument from \cite{Sibony}. Using Theorem~A,  $G_f^+$ is locally Hölder--continuous near $K^+(f)$; since $s(\CC) \cap K^+(f) = \Sigma_\kappa$ is a compact set, that property is global near the almost--sure spectrum and so $\gamma_\kappa$ is Hölder--continuous near $\Sigma_\kappa$. Denote by $\tau$ the exponent of H\"older continuity.

To show that $k_\kappa$ is H\"older continuous, consider two real numbers $E_2>E_1$. Let $M$ be the middle point of the segment $[E_1,E_2]$ and $R=\vert E_2-E_1\vert/2$ be the distance from $M$ to $E_1$. Denote by $D(r)\subset \CC$ the disk of radius $r$ centred at $M$. Let $\psi\colon \CC\to \RR_+$ be a smooth function which is equal to $1$ on $D(R)$ and equal to $0$ on $\CC\setminus D(2R)$, 
and whose partial derivatives of order $1$ and $2$ are bounded from above by $100 R^{-2}$ (such a function exists, see \cite{HarPol}). Then, 
\begin{eqnarray*}
\vert k_\kappa(E_2)-k_\kappa(E_1)\vert & = & \int_{[E_2,E_1]} \rmd k_\kappa(E) \\
& \leq & \int_{D(R)} \rmd \rmd ^c (\gamma_\kappa - \gamma_\kappa(M))\\
& \leq & \int_{D(3R)} \psi \cdot \rmd \rmd ^c(\gamma_\kappa - \gamma_\kappa(M))\\
& \leq & \int_{D(3R)} \rmd \rmd ^c\psi\cdot  (\gamma_\kappa - \gamma_\kappa(M))\\
& \leq & C^{st} R^\tau \mathrm{Area}(D(3R)) R^{-2}\\
& \leq & 9\pi C^{st} \vert E_2-E_1\vert ^\tau 
\end{eqnarray*}
for some uniform constant $C^{st}$   because $\gamma_\kappa$ is H\"older continuous (with exponent $\tau$) on a neighbourhood of $\Sigma_\kappa$. 

The same proof shows that $\rmd k_\kappa$ does not charge any closed subset of $\CC$ whose Hausdorff dimension is less than $\tau$
(see \cite{Sibony}).  \qed

\subsubsection{Hausdorff dimension of the density of states}\label{par:hausdorff-ds}
Once we know that $\gamma_\kappa$ is equal to $(\alpha+\beta)^{-1}G_f^+\circ s$, we can generalize the first results of Damanik and Gorodetski concerning the Hausdorff dimension of the density of states (proved in \cite{Dam4} for the Fibonacci substitution). Doing this, we obtain an alternative (but almost equal) proof of some of the results of May Mei (see \cite{Mei}). But first recall the following definition: we say that a finite measure $\mu$ on $\RR$ is of \textbf{exact dimension} $c \in \RR$ if for $\mu$--almost every $x \in \RR$ we have : 
\begin{displaymath}
\lim_{\varepsilon \to 0} \dfrac{\log \mu([x-\varepsilon,x+\varepsilon])}{\log\varepsilon} = c
\end{displaymath}
(see also \cite[p. 174]{Fal}).

\begin{theo}[Damanik, Gorodetski, Mei]
Let $\varphi$ be a positive and hyperbolic automorphism of the free group ${\mathbf{F}}_2$. Let $H_{\kappa, w}$ be the corresponding
family of discrete Schr\"odinger operators. 
For small coupling factors $0<\kappa< \kappa_0$, the density of states $\rmd k_\kappa$ is of exact dimension $\dim(\kappa)$, i.e. 
for $\rmd k_\kappa$-almost every real number $E$, 
\[
\lim_{\varepsilon\to 0} \dfrac{\log \rmd k_\kappa[E-\varepsilon, E+\varepsilon]}{\vert \log(\varepsilon)\vert} = \dim(\kappa).
\]
Moreover, 
\begin{itemize}
\item[(1)] $\dim(\kappa)$ is a ${\mathcal{C}}^\infty$-smooth function of $\kappa\in (0,\kappa_0)$;
\item[(2)] $\lim_{\kappa\to 0}\dim(\kappa)=1$;
\item[(3)] $\dim(\kappa) < \mathrm{H}_{\mathrm{dim}}(\Sigma_\kappa) < 1$ for $\kappa \in (0,\kappa_0)$, where $\mathrm{H}_{\mathrm{dim}}$ denotes the Hausdorff dimension of a given set;
\item[(4)] $\dim(\kappa)$ coincides with the infimum of $\mathrm{H}_{\mathrm{dim}}(S)$ 
over all measurable sets $S$ such that $\rmd k_\kappa(S)=1$. 
\end{itemize}
\end{theo} 

The proof is due to Damanik, Gorodetski, and Mei. Let us explain how one can relate its proof to Theorem~A and Theorem~B:\\

\noindent{\bf{a.--}} The dynamics of $f$ on the intersection of its filled Julia sets $K^+(f)\cap K^+(f^{-1})$ is {\sl{uniformly hyperbolic}}, 
the filled Julia set $K^+(f)$ is the support of a lamination by holomorphic curves, and the current $T^+_f$ is a current of integration on this lamination with respect to a transverse measure $\mu^+_f$ (see \cite{BHPS}).\\

\noindent{\bf{b.--}} The Schr\"odinger curve $s$ is {\sl{transverse to the lamination}} of $K^+(f)$ if the coupling factor is sufficiently small. This is proved in \cite{Dam5}; it follows from the transversality for $\kappa=0$ and a study of the bifurcation from $\kappa=0$ to 
$\kappa>0$.\\

\noindent{\bf{c.--}} There exists $\kappa_0'$ such that, for $0< \kappa < \kappa_0'$, {\sl{there are two saddle periodic points $p(\kappa)$ and $q(\kappa)$ of $f$ on $\Sr_{4+\kappa^2}$ with distinct multipliers}}. 

To prove this, take a periodic point $p$ on $\Sr_4$
which is not a singular point of $\Sr_4$. Deform it   into a family of periodic points $p(\kappa)$ for $-\kappa_1(p)\leq \kappa \leq \kappa_1(p)$. Do the same for a second periodic point $q$: it can be deformed into $q(\kappa)$ for $\kappa_1(q) < q < \kappa_1(q).$ If the multipliers of $p(\kappa)$ and $q(\kappa)$ are equal for a sequence of parameters $\kappa_n>0$ converging to $0$, they are equal for all $\kappa$ because they are analytic functions of $\kappa$. In particular, $q$ can be analytically deformed along the interval $[-\kappa_1(p),0]$. Thus, if the assertion was not satisfied, there would exist $\kappa_1>0$ such that all periodic points
of $f$ on $\Sr_4$ (distinct from the singularities) could be analytically deformed to saddle periodic points of the same period for $\kappa_1(p)<\kappa< 0$. This would contradict the fact that the topological entropy of $f$ on $\Sr_{4-\varepsilon}(\RR)$ is strictly less
than $\log(\lambda)$ for $\varepsilon >0$, a property that implies that most periodic points of $f$ on $\Sr_{4-\varepsilon}(\CC)$ are not real (see \cite{BHPS}).\\

With these three remarks in hand, one can then copy the proof given by Damanik and Gorodetski in \cite{Dam4}.

\subsubsection{Convergence theorems}\label{par:consequence-AS}
From \cite{BHPS} and \cite{BS1} (see also \cite{Sibony}, \cite{BS2}) one gets the following convergence theorem. {\sl{Let $f$ be a hyperbolic automorphism
of the surface $\Sr_D$. Let $T$ be a positive current and $\psi$ a smooth non-negative function with compact support which vanishes in a neighbourhood of the support of $\partial T$. Then, the sequence of currents 
\[
\dfrac{1}{\lambda^n}(f^n)^*(\psi T)
\]
converges towards a multiple $cT_f^+$, with $c=\langle T_f^-\vert \psi T\rangle$.}} For instance, $T$ can be the current of 
integration on an algebraic curve $C\subset \Sr_D$. 

Our goal is to explain, heuristically, why this result is similar to Avron-Simon convergence theorem for the density of states (see Theorem~\ref{thm:Avron-Simon}). 

Consider the restriction  $H_{\kappa,w}^N$ of the Schr\"odinger operator to some interval $[0,N]\subset \ZZ$. 
If $(u(0), \ldots, u(N))$ is an eigenfunction of $H_{\kappa,w}^N$ with eigenvalue $E$, then $(u(2),u(1))$ is obtained from $(u(1),u(0))$ by 
the linear action of the matrix $M_{\kappa}^E(w(0))$, \ldots, and $(u(N),u(N-1))$ is obtained from $(u(1),u(0))$ by 
the action of the product $M_{\kappa}^E(w(N-2)) \dots M_{\kappa}^E(w(0))$. 

Now, restrict the study to $w=u_+$, the infinite $\iota$-invariant word, and to intervals $[0,\ell(n)]$, where $\ell(n)$ is the length of the word $\iota^n(a)$. When $n$ goes 
to infinity, $\ell(n)$ behaves approximately like $\lambda^n$. 
With such a choice, the trace of the product $M_{\kappa}^E(w(\ell(n)-2)) \dots M_{\kappa}^E(w(0))$ is equal to the first coordinate of 
$f^n(s(E))$. Thus, if
\[
\Lambda:=\{E\; \vert\; \tr(M_{\kappa}^E(\iota^n(a)))=2\} ,
\]
then
\[
\Lambda=\{E\; \vert \; s(E)\in (f^n)^{-1}(C_2)\}
\]
where $C_2$ is the algebraic curve $C_2=\{(x,y,z)\in \Sr_D\vert x=2\}$. In other words, $\Lambda$ corresponds to the intersection of the algebraic
curve $s(\CC)$ with the algebraic curve $f^{-n}(C_2)$; it contains approximately $\lambda^n$ points, and the convergence
theorem for currents tells us, roughly, that the average measure on these $\lambda^n$ points converges towards $s^*(T^+_f)$, up to some multiplicative factor. 

On the other hand, the trace of a matrix $M\in SL(2,\RR)$ is $2$ if and only if $1$ is an eigenvalue of $M$. Thus, a complex number $E$ is in $\Lambda$ if and only if there is an eigenvector $(u(0), \ldots, u(\ell(n)))$ of $H_{\kappa,w}^{\ell(n)}$
with eigenvalue $E$ such that 
\[
(\star) \quad \left(\begin{array}{c} u(\ell(n)) \\ u(\ell(n)-1) \end{array}\right)=\left(\begin{array}{c} u(\ell(1)) \\ u(0) \end{array}\right);
\]
these are mixed boundary conditions (not the usual Dirichlet conditions, as in \cite{AS}). Thus, the convergence theorem for currents implies a convergence theorem for the density of states of $H^N_{\kappa,w}$ with the boundary conditions $(\star)$. Changing the
curve $C_2$ into another algebraic curve (for instance ${x=3}$), one gets different boundary conditions.

To sum up, Avron-Simon convergence theorem corresponds to the convergence theorem towards $T^+_f$, with the following differences: One only gets convergence along subsequences (one has to take $N=\ell(n)$), the boundary conditions
are not the classical ones, but one gets convergence theorems which are valid in $\Sr_D$ (not only along the Schr\"odinger curve) and work for all positive currents.

\ack

This paper was first drafted as a Masters thesis under the supervision of Serge Cantat and Frank Loray at the University of Rennes 1; the author would like to thank them both for their invaluable advice and flawless supervision. Paragraphs \ref{par:hausdorff-ds} and \ref{par:consequence-AS} were written by S. Cantat. The author would also like to acknowledge Christophe Dupont for his insight, Andrew Sale for his input regarding the English language, and Anton Gorodetski for his guidance. 

After a first draft of this paper circulated, A. Gorodetski provided the reference \cite{Mei} and William Yessen very kindly pointed out \cite{Y1,Y2} as having interesting ties to the author's work.

Funding for the author originated from the École Normale Supérieure de Cachan and the Université de Rennes 1.

\nocite{*}
\bibliographystyle{unsrt}
\bibliography{2013-green.bib}

\end{document}